\documentclass[11pt]{article}
\usepackage[utf8]{inputenc}
\usepackage[T1]{fontenc}
\usepackage{xspace}
\usepackage{amsthm}
\usepackage{amsmath}
\usepackage{amssymb}
\usepackage{bbm}
\usepackage{tikz}
\usepackage{subfigure}
\usepackage{keyval}
\usepackage[margin=3cm]{geometry}

\tikzstyle{nodino}=[circle,draw,fill,inner sep=0pt,minimum size=0.5mm]
\tikzstyle{infinito}=[circle,inner sep=0pt,minimum size=0mm]
\tikzstyle{nodo}=[circle,draw,fill,inner sep=0pt, minimum size=0.5*width("k")]

\newcommand{\rr}{{\mathbb R}}

\newcommand{\nn}{{\mathbb N}}
\newcommand{\G}{{\mathcal{G}}}
\newcommand{\udot}{\|u'\|_{L^2(\mathcal{G})}}
\newcommand{\uLp}{\|u\|_{L^p(\mathcal{G})}}
\newcommand{\uLtwo}{\|u\|_{L^2(\mathcal{G})}}
\newcommand{\HmuG}{H_\mu^1(\mathcal{G})}
\newcommand{\uLsix}{\|u\|_{L^6(\mathcal{G})}}
\newcommand{\uHone}{\|u\|_{H^1(\G)}}

\usetikzlibrary{graphs}

\theoremstyle{plain} 
\newtheorem{thm}{Theorem}[section]

\newtheorem{prop}[thm]{Proposition} 

\theoremstyle{definition}

\theoremstyle{definition} 
\newtheorem{defn}{Definition}[section] 

\theoremstyle{remark} 
\newtheorem{remark}{Remark}[section] 

\tikzset{every loop/.style={min distance=10mm,in=300,out=240,looseness=10}}

\tikzset{place/.style={circle,thick,draw=blue!75,fill=blue!20,minimum
size=6mm}}
\tikzset{place2/.style={circle,thick,draw=red!75,fill=red!20,minimum
size=6mm}}

\title{Existence of infinite stationary solutions of the $L^2$-subcritical and critical NLSE on compact metric graphs}

\author{Simone Dovetta \\{\small  Dipartimento di Scienze
		Matematiche ``G.L. Lagrange'', Politecnico di Torino } \\ {\small
		Corso Duca degli Abruzzi, 24, 10129 Torino, Italy}}

\begin{document}

\maketitle

\begin{abstract}
	We investigate the existence of stationary solutions for the Nonlinear Schr\"odinger equation on compact metric graphs. In the $L^2$-subcritical setting, we prove the existence of an infinite number of such solutions, for every value of the mass. In the critical regime, this infinity of solutions is established to exists if and only if the mass is lower or equal to a threshold value. Moreover, the relation between this threshold and the topology of the graph is characterized.
	The investigation is based on variational techniques and some new versions of Gagliardo-Nirenberg inequalities. 
\end{abstract}
	
	\section{Introduction}
		
		In this paper we discuss the existence of \textit{stationary solutions} for the NLS equation on a general \textit{compact metric graph} $\G$. In particular, we prove existence of critical points for the NLS energy functional
		\begin{equation}
			\label{NLS energy def}
			E(u,\G)=\frac{1}{2}\udot^2-\frac{1}{p}\uLp^p=\frac{1}{2}\int_{\G}|u'|^2dx-\frac{1}{p}\int_{\G}|u|^pdx
		\end{equation}
        \par\noindent with $p\in(2,6]$, under the mass constraint
        \begin{equation}
        	\label{mass constraint intro}
        	\uLtwo^2=\mu>0\quad.
        \end{equation}
        
        Such critical points solve, for suitable $\lambda\in\rr$, the stationary Schr\"odinger equation with the focusing pure power nonlinearity
        \begin{equation}
        \label{stationary NLS  intro}
        u''+|u|^{p-2}u=\lambda u
        \end{equation}
        \par\noindent on every edge of $\G$, with Kirchhoff conditions (see equation \eqref{Kirchhoff condition}) at the nodes. Throughout all this work, we limit ourselves to deal with real-valued functions.
        
        Our search for critical points of \eqref{NLS energy def} is twofold. On one hand, we investigate the existence of global minimizers of the constrained energy, called \textit{ground states}. Secondly, we turn our attention to a more general class of critical points, that do not need to be minimizers, usually called \textit{bound states}.
        
        We analyse both the subcritical regime $p\in(2,6)$ and critical one $p=6$, proving that the situation changes significantly. 
        
        In the subcritical case, we prove the existence of an infinite number of stationary solutions, with energy increasing to infinity, for every value of the mass and regardless of the topology of $\G$. This is stated in the following theorem.
        
        \begin{thm}
        	\label{subcritical compact}
        	Let $\G$ be a compact graph and $p\in(2,6)$. Then, for every $\mu>0$ there exist a ground state of \eqref{NLS energy functional} of mass $\mu$, and a sequence of bound states $\{u_k\}_{k\in\nn}$ of mass $\mu$, so that:
        	\begin{equation}
        		\label{energy to infty bound states subcritical}
        		E(u_k,\G)\to\infty\quad\quad\textit{for }k\to\infty.
        	\end{equation}
        \end{thm}
        
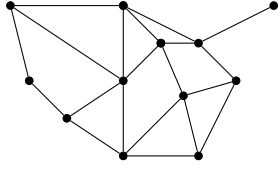
\begin{figure}[t]
		    	\centering
		    	\begin{tikzpicture}[xscale= 0.5,yscale=0.5] 
		    	\node at (-.5,2) [nodo] (02) {};
		    	\node at (2,2) [nodo] (22) {};
		    	\node at (2,4) [nodo] (24) {};
		    	\node at (3.6,1.6) [nodo] (42) {};
		    	\node at (3,3) [nodo] (33) {};
		    	\node at (5,2) [nodo] (52) {};
		    	\node at (3,3) [nodo] (32) {};
		    	\node at (4,3) [nodo] (43) {};
		    	\node at (-1,4) [nodo] (04) {};
		    	\node at (2,4) [nodo] (24) {};
		    	\node at (.5,1) [nodo] (11) {};
		    	\node at (2,0) [nodo] (20) {};
		    	\node at (4,0) [nodo] (40) {};
		    	
		    	\node at (6,4) [nodo] (term){};
		    	\draw [-] (43)--(term);
		    	
		    	\draw[-] (02)--(04);
		    	\draw[-] (04)--(24);
		    	\draw[-] (04)--(22);
		    	\draw[-] (24)--(22);
		    	\draw[-] (02)--(11);
		    	\draw[-] (11)--(22);
		    	\draw[-] (11)--(20);
		    	\draw[-] (20)--(22);
		    	\draw[-] (22)--(33);
		    	\draw[-] (24)--(33);
		    	\draw[-] (24)--(43);
		    	\draw[-] (33)--(43);
		    	\draw[-] (43)--(52);
		    	\draw[-] (33)--(42);
		    	\draw[-] (20)--(42);
		    	\draw[-] (20)--(40);
		    	\draw[-] (40)--(42);
		    	\draw[-] (42)--(52);
		    	\draw[-] (40)--(52);
		    	\end{tikzpicture}
		    	\caption{\small A compact graph with a terminal edge.}
		    	\label{grafo con terminal edge}
		    \end{figure}

        The critical regime is more subtle. We recall that by a \textit{terminal edge} (i.e. a tip, see Fig. \ref{grafo con terminal edge}) we mean an edge of $\G$ such that one of its endpoint is a vertex of degree one. It turns out that both ground state and bound states exist only for masses smaller than a threshold value that depends on whether the graph has at least one terminal edge or not. When $\G$ has a terminal edge, the threshold value of the mass is equal to the $L^2$-critical mass on the half-line $\rr^+$, $\mu_{\rr^+}=\sqrt{3}\pi/4$; in all other cases, it coincides with the $L^2$-critical mass on the whole real line $\rr$, $\mu_\rr=\sqrt{3}\pi/2$ (later on in this section we will briefly recall from where these two quantities arise).

        The theorem below establishes our main result at the critical exponent. 
        
        \begin{thm}
        	\label{critical thm compact}
        	Let $\G$ be a compact graph and $p=6$. Then:
        	\begin{itemize}
        		\item[(i)] if $\G$ has a terminal edge (Fig. \ref{grafo con terminal edge}), then a ground state of mass $\mu$ exists if and only if $\mu\leq\mu_{\rr^+}$. Moreover, for every $\mu<\mu_{\rr^+}$ there exists a sequence of bound states $\{u_k\}_{k\in\nn}$ of mass $\mu$ and $E(u_k,\G)\to\infty$;
        		\item[(ii)] if $\G$ has no terminal edge, then ground states of mass $\mu$ exist if and only if $\mu\leq\mu_{\rr}$. Moreover, for every $\mu<\mu_{\rr}$ there exists a sequence of bound states $\{u_k\}_{k\in\nn}$ of mass $\mu$ and $E(u_k,\G)\to\infty$.
        	\end{itemize}
        \end{thm}

        Roughly speaking, one may interpret our results by saying that the presence of a terminal edge forces $\G$ to share a critical behaviour similar to the half-line $\rr^+$, while in all other cases compact graphs seem to fake the line $\rr$. 
		
		The problem of the existence of a ground state on the real line $\G=\rr$ (that can be seen as two half-lines glued together at a single vertex) is nowadays classical (see \cite{cazenave}). In the subcritical regime, for every value of the mass $\mu$, ground states are unique up to translations and change of sign and are the so-called \textit{solitons} $\phi_\mu$, with strictly negative energy level. When $p=6$, on the contrary, the infimum of \eqref{NLS energy def} on $\rr$ undertakes a sharp transition from 0 to $-\infty$ when the value of $\mu$ exceeds $\mu_\rr$. A whole family $\{\phi_\lambda\}_{\lambda>0}$ of solitons realizing $E(\phi_\lambda,\rr)=0$ exists if and only if $\mu=\mu_{\rr}$.
		The same portrait holds when $\G=\rr^+$, with the threshold value of the mass becoming $\mu_{\rr^+}$ and the ground states are given by the \textit{half-solitons}, i.e. the restriction of $\phi_{2\mu}$ to $\rr^+$.

		 The study of nonlinear dynamics on quantum graphs has recently gained a considerable amount of interest, firstly motivated by a large class of physical applications, that range from optical networks to Bose-Einstein condensates. The evolution equation in this case is the time-dependent NLS
		 \begin{equation*}
		     i\partial_t\psi(t,x)=-\psi''(t,x)+g|\psi(t,x)|^{p-2}\psi(t,x)
		 \end{equation*}
		 \par\noindent that reduces to \eqref{stationary NLS  intro} when considering an attractive two-body interaction between the elementary components of the condensate (i.e. $g=-1$) and restricting to standing wave solutions $\psi(t,x)=e^{-i\lambda x}u(x)$ (for a more extended discussion of the physical interpretation see \cite{adamiserratilli_parma}).
		 
		 Former investigations were initiated in \cite{mehmeti,below,bona}, and after that, research has been pushed in several different directions. However, during the past years, the attention has been particularly focused on non-compact graphs with a finite number of nodes and at least a half-line. We point out that in such a physical context one should consider complex-valued functions; however, due to the invariance of eq. \eqref{stationary NLS  intro} under multiplication by a phase factor, one can restrict to the case of real-valued functions.

		 Evolution of solitary waves (\cite{squadraazzurra1}) and existence of ground states (\cite{squadraazzurra2,squadraazzurra3,squadraazzurra4}) for star graphs have both been addressed, while, for general graphs with half-lines, the investigations developed in \cite{adamiserratilli2014,adamiserratilli2015} for the subcritical case and in \cite{adamiserratilli2017_1} for the critical one revealed that both topological and metric properties of the domain definitely play a key role in allowing or preventing ground states from existing. Particularly, let us just mention that, at the critical exponent, for a certain class of graphs, ground states were proved to exist for a continuum of masses, similarly to what we established for compact graphs in Theorem \ref{critical thm compact}.
		 
		 Problems concerning bound states on graphs with half-lines have recently been addressed in \cite{adamiserratilli2017_2} considering minimization of the energy \eqref{NLS energy def} among functions of prescribed mass and fulfilling additional constraints.
		 
		 Moreover, existence of ground states and bound states on non-compact graphs was investigated also for the NLS equation with concentrated nonlinearity in \cite{tentarelli,serratentarelli2016_1,serratentarelli2016_2}. Specifically, the general scheme followed in \cite{serratentarelli2016_2} provides the tools we will use in this paper when dealing with bound states (see Section 2).
		 
		 Finally, we note that analysis of different classes of non-compact graphs has been recently initiated too, for instance in \cite{adamidovettaserratilli}, that deals with the ground state existence problem on an infinite, periodic graph (the two-dimensional grid).
		 
		 Considering compact graphs, the matter still seems to be quite unexplored. Something has been done on specific examples, for instance in \cite{Marzuola} and \cite{gustafson} for the cubic NLSE, but the techniques used there cannot be recovered to deal with the general exponent. Our contribution here consists in extending the scope of the analysis to more general graphs and nonlinearity powers. Particularly, we provide here a full topological characterization of the existence of ground states.

		 	Let us briefly recall some standard definitions and facts on compact metric graphs (for more details we refer for instance to \cite{adamiserratilli2014, berkolaiko, Kuchment}). Throughout the paper, a \textit{connected metric graph} $\G=(V,E)$ denotes a connected metric space built up by closed line intervals, the \textit{edges}, glued together with the identification of some of their endpoints, the \textit{nodes}. The peculiar way in which these identifications are performed defines the topology of the graph $\G$. Moreover, both multiple edges and self-loops are allowed. Every edge $e\in E$ is identified with an interval $I_e=[0,\ell_e]$, and a coordinate $x_e$ is chosen on $I_e$ providing orientation.
		 	
		 	A metric graph $\G$ is said to be \textit{compact} if and only if it has a finite number of nodes and edges and it has no edge of infinite length.
		 	
		 	Functions on a general graph $\G$ can be defined via their restriction to the edges. Indeed, every $u:\G\rightarrow\rr$ can be interpreted as a family of functions $(u_e)_{e\in E}$, where $u_e:I_e\rightarrow\rr$ is the restriction of $u$ to the edge represented by $I_e$. Thus, it is straightforward to define the usual functional spaces $L^p(\G)$ as 
		 	\begin{equation*}
		 	L^p(\G)=\{u:\G\rightarrow\rr:u_e\in L^p(I_e), \forall e\in E\}
		 	\end{equation*}
		 	\par\noindent and $H^1(\G)$ as
		 	\begin{equation*}
		 	H^1(\G)=\{u:\G\rightarrow\rr\text{ continuous}:u_e\in H^1(I_e),\forall e\in E\}
		 	\end{equation*}
		 	\par\noindent
		 	Note that requiring the continuity of $u$ implies that it is continuous at vertices. 
		 	
		 	Furthermore, we introduce the space
		 	\begin{equation*}
		 	H_\mu^1(\G)=\{u\in H^1(\G):\uLtwo^2=\mu\}
		 	\end{equation*}	 
		    \par\noindent that embodies the mass constraint \eqref{mass constraint intro}.
		    
		    The paper is organized as follows. In Section 2 we recall some standard notions from the Calculus of Variation we need for our purposes. In Section 3 we prove a general compactness property that holds for every compact graphs and we completely deal with the subcritical regime, proving Theorem \ref{subcritical compact}. In Section 4 we derive new technical tools playing a key role in the critical case, and we use them in Section 5 to provide the proof of Theorem \ref{critical thm compact}.

%

	\section{Variational framework}

	 As anticipated in the Introduction, our aim is to prove the existence of stationary solutions for the nonlinear Schr\"odinger equation with a pure power nonlinearity, defined as follows.
	 \begin{defn}
	 	\label{stationary solution NLS}
	 	Let $\G=(V,E)$ be a connected metric graph. A function $u\in H_\mu^1(\G)$ is said to be a \textit{stationary solution of mass $\mu$ for the NLS equation on $\G$ with Kirchhoff conditions at the nodes} if:
	 	\begin{itemize}
	 		\item[(i)] there exists $\lambda\in\rr$ such that, for every $e\in E$
	 		\begin{equation}
	 		\label{NLS on edge e}
	 			u_e''+|u_e|^{p-2}u_e=\lambda u_e,
	 		\end{equation}
	 		\item[(ii)] for every vertex $v\in V$,
	 		\begin{equation}
	 			\label{Kirchhoff condition}
	 			\sum_{e\succ v}\frac{du_e}{dx_e}(v)=0
	 		\end{equation}
	 	\end{itemize}
	 \end{defn}
	
	The symbol $e\succ v$ means that the sum above is extended to all edges $e$ incident at $v$, and $\frac{du_e}{dx_e}(v)$ can be either $u_e'(0)$ or $-u_e'(\ell_e)$, depending on the fact that $e$ reaches $v$ in $x_e=0$ or $x_e=\ell_e$. 
	
	Solutions of the stationary NLS equation as in Definition \ref{stationary solution NLS} are usually called \textit{bound states}, and a simple variational argument (see \cite{adamiserratilli2014}) shows that they can be found as constrained critical points of the \textit{energy functional} $E:H^1(\G)\rightarrow\rr$ defined as 
	\begin{equation}
		\label{NLS energy functional}
	     E(u,\G):=\frac{1}{2}\udot^2-\frac{1}{p}\uLp^p=\frac{1}{2}\int_{\G}|u'|^2dx-\frac{1}{p}\int_{\G}|u|^pdx
	\end{equation} 
	\par\noindent with the mass condition 
	
	\begin{equation*}
		\uLtwo=\mu.
	\end{equation*}
	
	Among all bound states, of particular interest are those critical points that globally minimize the NLS energy \eqref{NLS energy functional}, the so-called \textit{ground states}. We introduce the shorthand notation  
	\begin{equation}
	\label{ground state energy level function}
		\mathcal{E}(\mu):=\inf_{u\in H_\mu^1(\G)}E(u,\G),\quad\quad\mu\geq0
	\end{equation}

	In order to investigate the existence of bound states that are not in general global minimizers, we need some known results from the Calculus of Variations. Let us thus recall the following definition.
	\begin{defn}
         \label{Palais-Smale definition}
         Let $c\in\rr$. A sequence $\{u_n\}_{n\in\nn}\subset H_\mu^1(\G)$ is called a \textit{Palais-Smale sequence} for $E$ at level $c$ if, as $n\to\infty$,
         \begin{itemize}
         	\item[(i)] $E(u_n,\G)\to c$,
         	\item[(ii)] $\|E'(u_n,\G)\|_{T_{u_n}H_\mu^1(\G)}\to 0$ 
         \end{itemize}
         One says that $E$ satisfies the \textit{Palais-Smale condition} at level $c$ (denoted by $(PS)_c$) if every Palais-Smale sequence at level $c$ admits a subsequence strongly convergent in $H_\mu^1(\G)$.
         $E$ is said to satisfy the \textit{Palais-Smale condition} $(PS)$ if $(PS)_c$ holds for every $c$ admitting a Palais-Smale sequence.
	\end{defn}
	
   The following theorem, that is a unified version of some results presented in Section 10.2 of \cite{ambrosettimalchiodi}, states that the Palais-Smale condition ensures the existence of bound states.
	\begin{thm}
		\label{infinite bound states thm}
		Let $\mu>0$ and $E\in C^1(\HmuG,\rr)$ be as in \eqref{NLS energy functional}. Suppose that the Palais-Smale condition (PS) holds. Then there exists a sequence of bound states $\{u_k\}_{k\in\nn}\subset\HmuG$, and 
		\begin{equation}
			E(u_k,\G)\to+\infty\quad\quad\textit{for }k\to\infty
		\end{equation}
	\end{thm} 
	
	\begin{remark}
		In general, the results of \cite{ambrosettimalchiodi} summarized in the above theorem only imply that
		\begin{equation*}
			E(u_k,\G)\to\sup_{u\in\HmuG} E(u,\G).
		\end{equation*}
		\par\noindent However, it is immediate to see that, for every metric graph $\G$
		\begin{equation*}
		\sup_{u\in\HmuG} E(u,\G)=+\infty
		\end{equation*}
		\par\noindent thus justifying our statement of Theorem \ref{infinite bound states thm}.
	\end{remark}

    \section{Compactness, Gagliardo-Nirenberg inequalities and the subcritical regime}

    As one may expect, dealing with compact graphs reflects in a gain of compactness in $\HmuG$. The following proposition makes this fact more precise, as it establishes that, for Palais-Smale sequences, weak limits turn into strong ones.
    
    \begin{prop}
    	\label{compactness PS sequences thm}
    	Let $\G$ be a compact graph and $\{u_n\}_{n\in\nn}\subset\HmuG$ a Palais-Smale sequence for \eqref{NLS energy functional} bounded in $\HmuG$.
    	Then there exists $u\in\HmuG$ such that, up to subsequences, $u_n\to u$ strongly in $H^1(\G)$.
    \end{prop}
    
    \par\noindent\textit{Proof. }Since $\{u_n\}_{n\in\nn}$ is bounded in $\HmuG$, then $u_n\rightharpoonup u$ weakly in $H^1(\G)$, for some $u\in H^1(\G)$. By Sobolev compact embeddings, this implies (possibly passing to a subsequence)
    \begin{equation}
    	u_n\to u\quad\text{strongly in }L^p(\G),\quad p\geq 1
    \end{equation}
    \par\noindent so that $u\in\HmuG$.
    
    For every $u\in H_\mu^1(\G)$, we define the quantity
    \begin{equation}
    \label{definition of lambda}
    \lambda=\lambda(u):=-\frac{1}{\mu}E'(u)u=\frac{1}{\mu}\Big(\int_{\G}|u|^pdx-\int_{\G}|u'|^2dx\Big)
    \end{equation}
    \par\noindent and the linear functional $J(u):H^1(\G)\rightarrow\rr$
    \begin{equation}
    \label{action functional}
    J(u)v:=\int_{\G}u'v'dx-\int_{\G}|u|^{p-2}uv\,dx+\lambda\int_{\G}uv\,dx
    \end{equation}
    
    As proved in \cite{serratentarelli2016_2} (Section 2), since $\{u_n\}_{n\in\nn}$ is a bounded Palais-Smale sequence, condition (ii) in Definition \ref{Palais-Smale definition} can be conveniently rewritten as $J(u_n)\to0$ in $H^{-1}(\G)$. Moreover, if $\lambda_n:=\lambda(u_n)$ is as in \eqref{definition of lambda}, then $\{\lambda_n\}_{n\in\nn}$ is bounded in $\rr$, and (up to subsequences)
    \begin{equation}
    	\lambda_n\to\bar{\lambda}
    \end{equation}
    \noindent for some $\lambda\in\rr$.
    
    Let us now define the operator $A(u):H^1(\G)\to\rr$
    \begin{equation*}
    	\label{definition of A}
    	A(u)v:=\int_{\G}u'v'dx+\bar{\lambda}\int_{\G}uv
    \end{equation*}   
    \noindent and note that $A(u)(u_n-u)\to0$ by weak convergence of $u_n$ to $u$.
    
    Then we have:
    \begin{equation*}
    	\begin{split}
    	o(1)=&(J(u_n)-A(u))(u_n-u)\\
    	=&\int_{\G}|u_n'-u'|^2dx+\int_{\G}|u_n|^{p-2}u_n(u_n-u)dx+\lambda_n\int_{\G}u_n(u_n-u)dx\\
    	&-\bar{\lambda}\int_{\G}u(u_n-u)dx\\
    	=&\int_{\G}|u_n'-u'|^2dx+o(1)
    	\end{split}
    \end{equation*}
    \noindent by strong convergence in $L^p(\G)$ of $u_n$ to $u$ and by convergence of $\lambda_n$ to $\bar{\lambda}$, and this concludes the proof.\hspace{\stretch{1}} $\Box$
    
    One of the key tools in the study of the NLS energy functional \eqref{NLS energy functional} is the Gagliardo-Nirenberg inequality
    \begin{equation}
    \label{compact GN}
    \uLp^p\leq K_\G\uLtwo^{\frac{p}{2}+1}\uHone^{\frac{p}{2}-1}
    \end{equation}
    \par\noindent holding for every $p\geq 2$, $u\in H^1(\G)$ and any compact graph $\G$. Here $K_\G$ denotes the optimal constant.
   
   	Combining inequality \eqref{compact GN} and Proposition \ref{compactness PS sequences thm}, we are able to deal with the subcritical regime, proving Theorem \ref{subcritical compact}.
   	On the other hand, it turns out that the standard compact Gagliardo-Nirenberg inequality is not enough to manage the critical case. Theorem \ref{critical thm compact} thus requires modified versions of Gagliardo-Nirenberg inequalities that will be derived in the next section.

   \textit{Proof of Theorem \ref{subcritical compact}:} Plugging \eqref{compact GN} into \eqref{NLS energy functional} and using $p\leq6$, we have
   \begin{equation}
       \label{lower bound energy compact subcritical}
       \begin{split}
       E(u,\G)&\geq\frac{1}{2}\udot^2-\frac{1}{p}K_\G\mu^{\frac{p}{4}+\frac{1}{2}}\uHone^{\frac{p}{2}-1}\\
       &\geq\frac{1}{2}\udot^2\Big(1-C_1\mu^{\frac{p}{4}+\frac{1}{2}}\udot^{\frac{p}{2}-3}\Big)-C_2\mu^{\frac{p}{2}}
       \end{split}
   \end{equation}
  \par\noindent (where $C_1,C_2>0$ are proper constants), and so it follows
   \begin{equation}
   \label{infimum energy finite subcritical}
   	\mathcal{E}_\G(\mu)>-\infty
   \end{equation}
   \noindent for every $\mu>0$. 
   
   Let us now check that the Palais-Smale condition (PS) holds, for every value of the mass $\mu$. Fix $\mu>0$ and suppose $\{u_n\}_{n\in\nn}\subset\HmuG$ is a Palais-Smale sequence at level $c$. Then, by \eqref{lower bound energy compact subcritical}:
   \begin{equation*}
       c+o(1)=E(u_n,\G)\geq\frac{1}{2}\|u_n'\|_{L^2(\G)}^2\Big(1-C_1\mu^{\frac{p}{4}+\frac{1}{2}}\|u_n'\|_{L^2(\G)}^{\frac{p}{2}-3}\Big)-C_2\mu^{\frac{p}{2}}
   \end{equation*}
   \noindent Thus $\{u_n\}_{n\in\nn}$ is bounded in $H^1(\G)$. Hence, by Proposition \ref{compactness PS sequences thm}, $u_n\to u$ strongly in $H^1(\G)\cap L^p(\G)$, for some $u\in\HmuG$. Since this is true for every $c$, (PS) is proved.
   
   Since it always exists a minimizing sequence that is also a Palais-Smale sequence at level $c=\mathcal{E}_\G(\mu)$ (that is finite due to \eqref{infimum energy finite subcritical}), ground states exist for every $\mu>0$. Moreover, Theorem \ref{infinite bound states thm} applies and a sequence of bound states $\{u_k\}_{k\in\nn}$ exists for every value of the mass and
   \begin{equation*}
   	E(u_k,\G)\to\infty\quad\quad\quad\text{for }k\to\infty.
   \end{equation*}\hspace{\stretch{1}} $\Box$
   
   \section{Modified Gagliardo-Nirenberg inequalities}
    
   In order to deal with the critical case $p=6$, we state here a modified version of Gagliardo-Nirenberg inequalities that holds for general compact graph. The importance of the inequality in this form is that the $H^1-$norm of the function $u$ appearing in \eqref{compact GN} is replaced by the $L^2$-norm of the first derivative $u'$, as in the non-compact case (see for instance \cite{adamiserratilli2015}). Moreover, the constant involved here instead of $K_\G$ directly relates the mass value $\mu$ to $\mu_{\rr^+}$ or $\mu_\rr$, depending on the possible presence of a terminal edge. 
   
   A similar modified Gagliardo-Nirenberg inequality was derived for the first time in Lemma 4.4 in \cite{adamiserratilli2017_1} for general non-compact graphs, so, to some extent, our proof is an adaptation of the one given in \cite{adamiserratilli2017_1}.
   
   First, let us recall the standard Gagliardo-Nirenberg inequality with $p=6$ on $\rr$ and $\rr^+$
   
   \begin{equation}
   	\label{GN on r and r+}
   	\uLsix^6\leq K\uLtwo^{4}\udot^{2}
   \end{equation}
   
   \par\noindent where $K=\frac{3}{\mu_\rr^2},\frac{3}{\mu_{\rr^+}^2}$ denotes the optimal constant on $\rr$ and $\rr^+$ respectively.

   \begin{prop}
   	\label{modified GN tip prop}
   	Assume $\G$ is a compact graph with at least one terminal edge. Fix $\mu\in(0,\mu_{\rr^+}]$ and let $u\in\HmuG$. Then there exists a number $\theta=\theta(u)\in[0,\mu]$ such that
   	\begin{equation}
   		\label{modified GN tip formula}
   		\uLsix^6\leq3\Big(\frac{\mu-\theta}{\mu_{\rr^+}}\Big)^2\udot^2+C\theta^{1/2}
   	\end{equation}
   	\noindent where $C>0$ is a constant that depends only on $\G$.
   \end{prop}

   \textit{Proof:} If $u$ is constant on $\G$, then the result is immediate, so let us consider a non-constant function $u$. Replacing $u$ with $|u|$, we may assume $u>0$.
   
   Let $\ell:=|\G|$ be the total length of $\G$. Considering the decreasing rearrangement $u^*$ of $u$ on $[0,\ell)$, it is well-known (see e.g. \cite{adamiserratilli2014}) that $u^*\in H_\mu^1(0,\ell)$, $u^*$ is non-increasing on $[0,\ell]$ and
   \begin{equation}
       \label{monotone rearrangement norms}
       \begin{split}
          \|u^*\|_{L^6(0,\ell)}&=\uLsix\\
          \|(u^{*})'\|_{L^2(0,\ell)}&\leq\udot
       \end{split}
   \end{equation}
   
   Rephrasing Step 2 of the proof of Lemma 4.4. in \cite{adamiserratilli2017_1}, one can construct a function $v\in H^1(\rr^+)$ so that, for some $\theta=\theta(u)\in[0,\mu]$:
   \begin{itemize}
   	\item[(i)] $v(0)=u^*(0)$;
   	\item[(ii)] $\int_{0}^{\infty}|v|^2dx=\int_{0}^{\ell}|u^*|^2dx-\theta=\mu-\theta$;
   	\item[(iii)] $\int_{0}^{\infty}|v'|^2dx\leq\int_{0}^{\ell}|(u^*)'|^2dx+C\theta^{1/2}$;
   	\item[(iv)] $\int_{0}^{\infty}|v|^6dx\geq\int_{0}^{\ell}|u^*|^6dx-C\theta$;
   \end{itemize}
   \noindent with the constant $C>0$ depending only on $\G$.
   
   Now, by Gagliardo-Nirenberg inequality \eqref{GN on r and r+} on $\rr^+$, 
   \begin{equation}
   	\label{GN rr+ on v}
   	\|v\|_{L^6(\rr^+)}^6\leq\frac{3}{\mu_{\rr^+}^2}\|v\|_{L^2(\rr^+)}^4\|v'\|_{L^2(\rr^+)}^2=3\Big(\frac{\mu-\theta}{\mu_{\rr^+}}\Big)^2\|v'\|_{L^2(\rr^+)}^2\quad.
   \end{equation}

   Using properties (iii)-(iv) of $v$, we can then get back to $u^*$, and then to $u$. Indeed, by (iii) and \eqref{monotone rearrangement norms} we have
   \begin{equation}
   	\label{v 6 with u 6}
   	\|v\|_{L^6(\rr^+)}^6\leq\uLsix^6-C\theta
   \end{equation}
   \noindent while (iv) and \eqref{monotone rearrangement norms} again give
   \begin{equation}
   	\label{v' 2 with u' 2}
   	\|v'\|_{L^2(\rr^+)}^2\leq\udot^2+C\theta^{1/2}
   \end{equation}
   
   Finally, plugging \eqref{v 6 with u 6} and \eqref{v' 2 with u' 2} into \eqref{GN rr+ on v} and properly changing $C$, we get \eqref{modified GN tip formula} and the proof is complete.\hspace{\stretch{1}} $\Box$
   
   The following proposition is the analogue of the previous one if there is no tip in $\G$, where $\mu_\rr$ takes the place of $\mu_{\rr^+}$.
   \begin{prop}
      \label{modified GN no tip prop}
      Assume $\G$ is a compact graph with no terminal edge. Fix $\mu\in(0,\mu_{\rr}]$ and let $u\in\HmuG$. Then there exists a number $\theta=\theta(u)\in[0,\mu]$ such that
      \begin{equation}
      	\label{modified GN no tip formula}
      	\uLsix^6\leq3\Big(\frac{\mu-\theta}{\mu_{\rr}}\Big)^2\udot^2+C\theta^{1/2}
      \end{equation}
      \noindent where $C>0$ is a constant that depends only on $\G$.	
   \end{prop}
   
   \begin{figure}
   	\centering	
   	\begin{tikzpicture}[xscale= 0.9,yscale=0.7]
   	\filldraw[black] (0,2) circle (2 pt);
   	\filldraw[black] (-3,2) circle (2 pt);
   	\filldraw[black] (2,2) circle (2 pt);
   	\filldraw[black] (5,1.5) circle (2 pt);
   	\filldraw[black] (-4,0.5) circle (2 pt);
   	\filldraw[black] (-3.6,-0.5) circle (2 pt);
   	\draw (-1.5,2) ellipse (42 pt and 13 pt);
   	\draw (0,2)--(2,2);
   	\draw [-] (2,2) to [out=-45,in=-110] (5,1.5);
   	\draw [-] (2,2) to [out=-5,in=140] (5,1.5);
   	\draw (-3,2)--(-4,0.5);
   	\draw (-4,0.5)--(-3.6,-0.5);
   	\draw (-3.6,-0.5)--(-3,2);
   	\draw [-] (0,5) to [out=-45,in=180] (0.8,2.5);
   	\draw [-] (0.8,2.5) to [out=0,in=180] (4,7);
   	\draw [-] (4,7) to [out=0,in=90] (5,5.7);
   	\draw [-] (5,5.7) to [out=-90,in=-45] (2,3.4);
   	\draw [-] (0,5) to [out=-60,in=-45] (-3,3);
   	\draw [-] (0,5) to [out=135,in=-20] (-3,3);
   	\draw [-] (-3,3) to [out=160,in=0] (-3.5,4);
   	\draw [-] (-3.5,4) to [out=180,in=160] (-4,3);
   	\draw [-] (-4,3) to [out=-20,in=-135] (-3,3);
   	\draw [dashed] (5,1.5)--(5,5.7);
   	\draw [dashed] (4,7)--(4,2);
   	\draw [dashed] (0.8,2.5)--(0.8,2);
   	\node at (3.8,1.20) [label= \text{\footnotesize$x_0$}] {};
   	\draw [dashed] (0.9,2.07)--(2,2.07);
   	\draw [dashed] (2,2.05) to [out=5,in=0](3.8,2.05);
   	\node at (2.5,1.9) [label= \text{\large$\Sigma$}] {};
   	\draw [dashed] (4.1,1.9) to [out=170,in=130](4.9,1.5);
   	\node at (4.45,0.7) [label=\text{\large$\Gamma$}]{};
   	\end{tikzpicture}
   	\caption{Construction of paths $\Sigma$ and $\Gamma$ in the proof of Proposition \ref{modified GN no tip prop}.}
   	\label{proof picture}
   \end{figure}
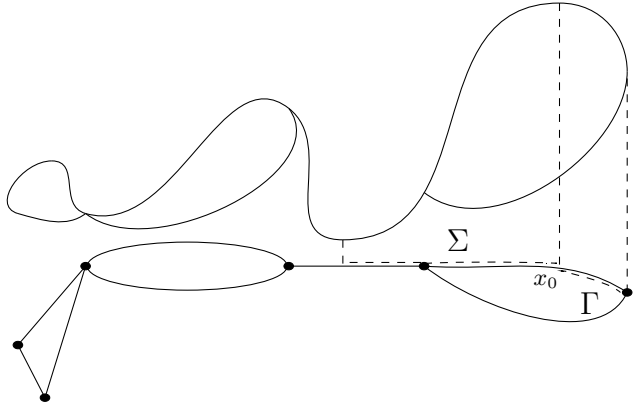

   \textit{Proof. } The line of the proof is similar to the one of Proposition \ref{modified GN tip prop}. The main difference is that it is possible to make use of the standard Gagliardo-Nirenberg inequality on $\rr$ instead of $\rr^+$.
   Again, we can consider $u$ not everywhere constant on $\G$, $u>0$ on $\G$. 
   
   
   Let us denote by $2\gamma$ the length of the shortest loop in $\G$. Suppose that $u$ realizes its maximum at some point $x_0$, and let $\Sigma$ be the shortest path in $\G$ going from $x_0$ to a point in which $u$ attains its minimum (it exists since $\G$ is connected). Since there is no terminal point in $\G$, $\Sigma$ can be extended for an additional length $\gamma$ from $x_0$ in the direction opposite to the minimum. We call $\Gamma$ the path starting at $x_0$ of length $\gamma$ obtained this way (see Figure \ref{proof picture}). Since $\G$ is connected and $u$ is continuous, it follows that $u_{\mid\G/\Gamma}$ reaches all the values in the range of $u$ at least once. Hence, rearranging monotonically $u_{\mid\Gamma}$ on $[0,\gamma]$ and $u_{\mid\G/\Gamma}$ on $(\gamma-\ell,0]$ and gluing them together at the origin, we get a function $\varphi\in H^1(\gamma-\ell,\gamma)$ such that:
   \begin{itemize}
   	\item[(a)] $\varphi(0)=\|\varphi\|_{L^\infty(\gamma-\ell,\gamma)}=\|u\|_{L^\infty(\G)}$ and $\varphi(\gamma-\ell)=m$;
   	\item[(b)] $\varphi$ is monotonically increasing on $[\gamma-\ell,0]$ and monotonically decreasing on $[0,\gamma]$;
   	\item[(c)] $\int_{\gamma-\ell}^{\ell}|\varphi|^2dx=\mu$;
   	\item[(d)] $\int_{\gamma-\ell}^{\ell}|\varphi|^6dx=\int_{\G}|u|^6dx$, \quad while $\int_{\gamma-\ell}^{\ell}|\varphi'|^2dx\leq\int_{\G}|u'|^2dx$.
   \end{itemize}
   
   Now, starting from $\varphi$, we construct a function $w\in H^1(\rr)$ and $\theta\in[0,\mu]$ satisfying:
   \begin{itemize}
   	\item[1.] $\int_{\rr}|w|^2dx=\mu-\theta$;
   	\item[2.] $\int_{\rr}|w'|^2dx\leq\int_{\gamma-\ell}^{\ell}|\varphi'|^2dx+C\theta^{1/2}$;
   	\item[3.] $\int_{\rr}|w|^6dx\geq\int_{\gamma-\ell}^{\ell}|\varphi|^6dx-C\theta$
   \end{itemize}
   \noindent where $C>0$ depends once again only on $\G$. Indeed, applying Step 2 of the proof of Lemma 4.4 in \cite{adamiserratilli2017_1} independently to $\varphi_{\mid[\gamma-\ell,0]}$ and $\varphi_{\mid[0,\gamma]}$, we get two functions $v_1,v_2\in H^1(\rr^+)$ satisfying properties (i)-(iv) with some non negative $\theta_1,\theta_2$ so that $\theta_1+\theta_2\in[0,\mu]$. Gluing together $v_1(x)$ and $v_2(-x)$ at $x=0$ and setting $\theta:=\theta_1+\theta_2$ we get $w$ as above. Here $\theta=\theta(u)$ depends by construction on the original function $u$.
   
   Applying Gagliardo-Nirenberg \eqref{GN on r and r+} with $K=\frac{3}{\mu_\rr^2}$ to $w$ and combining with 2.-3. and the properties of $\varphi$, we have
   \begin{equation*}
   	\uLsix^6-C\theta\leq\|w\|_{L^6(\rr)}^6\leq3\Big(\frac{\mu-\theta}{\mu_\rr}\Big)^2\|w'\|_{L^2(\rr)}^2\leq3\Big(\frac{\mu-\theta}{\mu_\rr}\Big)^2\big(\udot^2+C\theta^{1/2}\big)
   \end{equation*}
   \noindent Rearranging terms and modifying $C$ if necessary conclude the proof.\hspace{\stretch{1}} $\Box$

   \begin{remark}
   	Proposition \ref{modified GN no tip prop} states a stronger result than Proposition \ref{modified GN tip prop}, in the sense that it holds for a wider interval of masses, since $\mu_{\rr}>\mu_{\rr^+}$. This difference relies on the fact that, if $\G$ has no terminal edges, then every function on the graph shares at least two pre-images for almost every value it achieves. On the contrary, the presence of a terminal edge allows to exhibit functions on $\G$ with only one pre-image for an interval of values. Such a difference is strictly related to the properties of rearrangements and it is where it enters the field through our proofs. For a detailed overview on the connection with rearrangements we refer to \cite{adamiserratilli2015}.
   \end{remark}
   
   \section{Proof of the main result}
   
   We present here the proof of Theorem \ref{critical thm compact}, that takes advantage of the modified Gagliardo-Nirenberg inequalities derived before.
   
   \textit{Proof of Theorems \ref{critical thm compact}: } 
   Let us begin by assuming that $\G$ is a compact graph with at least one terminal edge.
   
   We first deal with the case $\mu>\mu_{\rr^+}$. Let $v\in H_{\mu}^1(\rr^+)$ so that $supp(v)=[0,1]$ and $E(v,\rr^+)<0$ (it surely exists since $E$ is unbounded from below on $\rr^+$ when $\mu>\mu_{\rr^+}$ \cite{cazenave}). Defining, for $\lambda>0$,
   \begin{equation}
   \label{u n to -infty}
   v_\lambda(x):=\sqrt{\lambda}v(\lambda x),\quad\quad x\in\rr^+
   \end{equation}
   \noindent we get a family of functions so that:
   \begin{equation*}
   \begin{split}
   \int_{\rr^+}|v_\lambda|^2dx&=\mu\\
   supp(v_\lambda)&=\Big[0,\frac{1}{\sqrt{\lambda}}\Big]\\
   E(v_\lambda,\rr^+)&=\lambda^2E(v,\rr^+)
   \end{split}
   \end{equation*}
   
   Thus, when $\lambda$ is large enough, the support of $v_\lambda$ can be contained by any edge of $\G$. Therefore, suppose $e$ is an edge of $\G$ and choose $\lambda>0$ so that $supp(v_\lambda)\subset I_e$. Then, defining $w_\lambda\in\HmuG$ as follows
   \begin{equation*}
   w_\lambda(x)=\begin{cases}
   v_\lambda(x) & \text{if }x\in I_e\\
   0 & \text{elsewhere on }\G
   \end{cases}
   \end{equation*}
   \noindent one gets $E(w_\lambda,\G)\to-\infty$ as $\lambda\to\infty$, thus
   \begin{equation}
   \label{unbounded inf mu greater mu r+}
   \mathcal{E}_\G(\mu)=-\infty
   \end{equation}
   \noindent for every $\mu>\mu_{\rr^+}$, and ground states do not exist.
   
   On the contrary, in order to ensure existence of ground states and bound states when $\mu\leq\mu_{\rr^+}$, it is sufficient to check the validity of the Palais-Smale condition. 
   
   Let us first consider $\mu<\mu_{\rr^+}$. Since $\theta\leq\mu$, for every $u\in\HmuG$, plugging inequality \eqref{modified GN tip formula} in the energy functional \eqref{NLS energy functional} one gets rid of $\theta$,
   \begin{equation}
   	\begin{split}
   	\label{lower bound critical GN tip}
   	E(u,\G)\geq\frac{1}{2}\udot^2\Big(1-\frac{\mu^2}{\mu_{\rr^+}^2}\Big)-C\mu^{1/2}
   	\end{split}
   \end{equation}
   \noindent thus implying
   \begin{equation}
   	\label{inf sup tip enegy critical}
   	\mathcal{E}_\G(\mu)>-\infty.
   \end{equation}
   
    Let $\{u_n\}_{n\in\nn}\subset\HmuG$ be a Palais-Smale sequence at level $c\in\rr$. Then, by $E(u_n)\to c$,
    \begin{equation}
      \label{PS critical tip}
      c+o(1)=E(u_n,\G)\geq\frac{1}{2}\|u_n\|^2\Big(1-\frac{\mu^2}{\mu_{\rr^+}^2}\Big)-C\mu^{1/2}
    \end{equation}
    \noindent and $\{u_n\}_{n\in\nn}$ is bounded in $H^1(\G)$. Then, by Proposition \ref{compactness PS sequences thm}, $u_n\to u$ strongly, for some $u\in\HmuG$. Since this holds for every $c$, (PS) is proved.

   Let us now consider the case $\mu=\mu_{\rr^+}$.  We denote here by $\theta_u$ the constant appearing in \eqref{modified GN tip formula} to stress its dependence on $u$. Even though it is no longer possible to ignore the role of $\theta$ as in \eqref{lower bound critical GN tip}, we are still able to deal with functions $u$ that realize strictly negative energy. Indeed, if $E(u,\G)=-\alpha<0$, then by \eqref{modified GN tip formula} one has
   \begin{equation}
   \label{lower bound mu rr+ critical}
   	\frac{1}{2}\udot^2\Big(1-\frac{\theta_u^2}{\mu_{\rr^+}^2}\Big)-C\theta_u^{1/2}\leq E(u,\G)=-\alpha<0
   \end{equation}
   \noindent showing that $\theta_u$ is bounded away from 0.
   
   On one hand, this ensures again that
   \begin{equation}
   	\label{infimum mu rr+ critical bounded}
   	\mathcal{E}_\G(\mu_{\rr^+})>-\infty
   \end{equation}
    
  On the other hand, it implies that the Palais-Smale condition can be recovered at least at negative leves $c<0$, thanks to the lower bound provided by inequality in \eqref{lower bound mu rr+ critical}. Therefore, if $\{u_n\}_{n\in\nn}\subset H_{\mu_{\rr^+}}^1(\G)$ is a Palais-Smale sequence at $c<0$, then $u_n\to u$ strongly in $H_{\mu_{\rr^+}}^1(\G)$. 
  
  Remember that it is not restrictive to consider minimizing sequences that are also Palais-Smale sequences at level $c=\mathcal{E}_\G(\mu_{\rr^+})$, and $c>-\infty$ for every $\mu\leq\mu_{\rr^+}$ by \eqref{inf sup tip enegy critical}-\eqref{infimum mu rr+ critical bounded}. 
  
  Moreover, the constant function $\varsigma\in\HmuG$ on $\G$
       
   \begin{equation}
      \label{constant on G}
      \varsigma:=\sqrt{\frac{\mu}{\ell}}
   \end{equation}

  \par\noindent ensures that, on every compact graph:
  \begin{equation*}
  	\mathcal{E}_\G(\mu_{\rr^+})\leq E(\varsigma,\G)=-\frac{\mu_{\rr^+}^3}{6\ell^2}<0
  \end{equation*}
  
  Since the Palais-Smale condition always holds at negative levels for $\mu\leq\mu_{\rr^+}$, if $\{u_n\}_{n\in\nn}\subset\HmuG$ is a minimizing Palais-Smale sequence, then it is compact in $\HmuG$, and ground states exist for every value of the mass less or equal to $\mu_{\rr^+}$. This proves the first part of (i) in Theorem \ref{critical thm compact}. 
  \noindent Moreover, since (PS) was proved to hold whenever $\mu$ is strictly less than $\mu_{\rr^+}$, then Theorem \ref{infinite bound states thm} applies in this case, and the second part of (i) in Theorem \ref{critical thm compact} follows.
  
  Statement (ii) follows the same argument, using the proper modified Gagliardo-Nirenberg \eqref{modified GN no tip formula} instead of \eqref{modified GN tip formula} whenever needed.\hspace{\stretch{1}} $\Box$
  
  \begin{remark}
  	Theorem \ref{critical thm compact} does not apply when $\mu=\mu_{\rr^+}$ and $\G$ has a terminal edge or $\mu=\mu_\rr$ and $\G$ has no terminal edge. The problem in this situation is that we do not know if the Palais-Smale condition holds at every level.
    Indeed, at the threshold mass, it is still true (thanks to the methods we used) that this condition holds at negative levels $c<0$.
  	Unfortunately, taking for instance sequences of functions with compact support approximating different half-solitons (resp. solitons), it is easy to see that it fails at level $c=0$ when $\mu=\mu_{\rr^+}$ (resp. $\mu=\mu_\rr$) on a graph with a terminal edge (resp. with no terminal edge). Moreover, when $c>0$, we are not able to establish whether Palais-Smale sequences are compact or may not converge.
  \end{remark}
  
   \begin{remark}
       Let $\ell:=|\G|$ denote the total length of $\G$. Using standard methods in the theory of stability developed in \cite{GSS}, it is possible to prove that, both in the subcritical and in the critical domain, there exists a threshold value of the mass, say $\mu^*=\mu^*(\G,p)$, such that the constant function $\varsigma\in\HmuG$ as in \eqref{constant on G} is a local minimum for \eqref{NLS energy functional} on the constrained manifold $\|u\|_{L^2(\G)}^2=\mu$, whenever $\mu<\mu^*$. Even though partial results in highly specific cases have been developed in \cite{Marzuola} and \cite{gustafson}, understanding when $\varsigma$ is actually a ground state of the energy on a general compact graph seems to lay a quite involved question.
    \end{remark}
  
  \medskip
  \medskip
  
  \par\noindent\textbf{Acknowledgements}
  \medskip
  
  \par\noindent The author is grateful to Enrico Serra for all enlightening suggestions he provided during the preparation of this work and to Riccardo Adami for constant support.


\begin{thebibliography}{99}
   	
   	  \bibitem{squadraazzurra1}
   	  Adami R., Cacciapuoti C., Finco D., Noja D.,
   	  \textit{Fast solitons on star graphs},
   	  Rev. Math. Phys. \textbf{23}(4), 409-451 (2011).
   	  
   	  \bibitem{squadraazzurra2}
   	  Adami R., Cacciapuoti C., Finco D., Noja D.,
   	  \textit{On the structure of critical energy levels for the cubic focusing NLS on star graphs},
   	  J. Phys. A \textbf{45}(19), 192001 (2012).
      
      \bibitem{squadraazzurra3}
   	  Adami R., Cacciapuoti C., Finco D., Noja D.,
   	  \textit{Variational properties and orbital stability of standing waves for NLS equation on a star graph},
   	  J. Diff. Eq. \textbf{257}(10), 3738-3777 (2014).
   	  
   	  \bibitem{squadraazzurra4}
   	  Adami R., Cacciapuoti C., Finco D., Noja D.,
   	  \textit{Stable standing waves for a NLS on star graphs as local minimizers of the constrained energy},
   	  J. Diff. Eq. \textbf{260} (10), 7397-7415 (2016).
      
      \bibitem{adamiserratilli2014}
      Adami  R., Serra E., Tilli P.,
      \textit{NLS ground states on graphs},
      Calc. Var. and PDEs \textbf{54(1)} (2015), 743–761.
      
      \bibitem{adamiserratilli2015}
      Adami R., Serra E., Tilli P.,
      \textit{Threshold phenomena and existence results for NLS ground states on graphs},
      Journal of Functional Analysis \textbf{271(1)} (2015).
      
      \bibitem{adamiserratilli2017_1}
      Adami R., Serra E., Tilli P.,
      \textit{Negative energy ground states for the $L^2$-critical NLSE on metric graphs},
      Commun. Math. Phys. \textbf{352} (2017), no. 1, 387-406.
      
      \bibitem{adamiserratilli2017_2}
      Adami R., Serra E., Tilli P.,
      \textit{Multiple positive bound states for the subcritical NLS equation on metric graphs},
      arXiv: 1706.07654.
      
      \bibitem{adamiserratilli_parma}
      Adami R., Serra E., Tilli P.,
      \textit{Nonlinear dynamics on branched structures and networks},
      Riv. Mat. Univ. Parma, Vol. 8, No. 1, 109-159 (2017)     
            
      \bibitem{adamidovettaserratilli}
      Adami R., Dovetta S., Serra E., Tilli P.,
      \textit{NLS ground states on the two-dimensional grid: dimensional crossover and a continuum of critical exponents} 
      (to appear)
      
      
      \bibitem{mehmeti}
      Ali Mehmeti F.,
      \textit{Nonlinear wavees in networks}, 
      Akademie Verlag Berling (1994).
      
      \bibitem{ambrosettimalchiodi}
      Ambrosetti A., Malchiodi A.,
      \textit{Nonlinear analysis and semilinear elliptic problems},
      Cambridge Studies in Advanced Mathematics, 104, Cambridge University Press, Cambridge, (2007).
      
      \bibitem{below}
      Below J. von, 
      \textit{An existence result for semilinear parabolic network equations with dynamical node conditions}, 
      Pitman Research Notes in Mathematical Series 266, 274-283, Longman, Harlow Essex (1992).
      
      \bibitem{berkolaiko}
      Berkolaiko G., Kuchment P.,
      \textit{Introduction to quantum graphs},
      Mathematical Surveys and Monographs, 186. AMS, Providence, RI (2013)
      
      \bibitem{bona}
      Bona J., Cascaval R.C., 
      \textit{Nonlinear dispersive waves on trees},
      Can. J. App. Math. \textbf{16}, 1-18 (2008).
      
     \bibitem{cazenave}
      Cazenave T.,
      \textit{Semilinear Schr\"odinger Equations},
      Courant Lecture Notes 10. American Mathematical Society, Providence, RI (2003).
      
      \bibitem{GSS}
      Grillakis M., Shatah J., Strauss W.,
      \textit{Stability theory of solitary waves in the presence of symmetry, I}, 
      J. Funct. An. \textbf{74} (1987), no. 1, 160-197.
      
      \bibitem{gustafson}
      Gustafson S., Le Coz S., Tsai T.P.,
      \textit{Stability of periodic waves of 1D cubic nonlinear Schr\"odinger equations},
      arXiv:1606.04215.
      
      \bibitem{Kuchment}
      Kuchment P.,
      \textit{Quantum graphs. I. Some basic structures},
      Waves Random Media \textbf{14} (2004), no. 1, 107-128.
      
      \bibitem{Marzuola}
      Marzuola J., Pelinovsky D.E.,
      \textit{Ground states on the dumbbell graph},
      arXiv: 1509.04721
      
      
      \bibitem{serratentarelli2016_2}
      Serra E., Tentarelli L.,
      \textit{Bound states of the NLS equation on metric graphs with localized nonlinearities},
      J. Diff. Eq. \textbf{260} (2016), no. 7, 5627-5644.
      
      \bibitem{serratentarelli2016_1}
      Serra E., Tentarelli L.,
      \textit{On the lack of bound states for certain NLS equations on metric graphs},
      Nonlinear Anal. \textbf{145} (2016), 68-82.
      
      \bibitem{tentarelli}
      Tentarelli L.,
      \textit{NLS ground states on metric graphs with localized nonlinearities},
      J. Math. Anal. Appl. \textbf{433} (2016), no. 1, 291-304.
      
      
   \end{thebibliography}
\end{document}